\documentclass{amsart}
\usepackage{amsmath, amsthm, amssymb}
\usepackage{pdfsync}
\input xy
\xyoption{all}

\theoremstyle{plain} 
\newtheorem{theorem} {Theorem}
\newtheorem{proposition}[theorem] {Proposition}
\newtheorem{lemma}[theorem] {Lemma}
\newtheorem{corollary}[theorem] {Corollary}

\theoremstyle{definition}

\newtheorem{definition}[theorem] {Definition}

\theoremstyle{remark}
\newtheorem{remark}[theorem] {Remark}

\newcommand{\mpair}[1]{\pair{\,#1\,}}

\newcommand{\mset}[1]{\set{\,#1\,}}

\newcommand{\pair}[1]{\langle #1\rangle}

\newcommand{\set}[1]{\{#1\}}

\newcommand{\rf}{root function\ }
\newcommand{\rfs}{root functions\ }
\newcommand{\rfp}{root function.}

\newcommand{\inv}{^{-1}}

\newcommand{\wh}{\widehat}

\newcommand{\nc}{\newcommand}

\nc{\on}{\operatorname}

\nc{\Z}{{\mathbb Z}}
\nc{\C}{{\mathbb C}}
\nc{\R}{{\mathbb R}}
\nc{\bbP}{{\mathbb P}}
\nc{\bF}{{\mathbb F}}

\nc{\boldD}{{\mathbb D}}
\nc{\oo}{{\mf O}}
\nc{\N}{{\mathbb N}}
\nc{\bib}{\bibitem}
\nc{\pa}{\partial}
\nc{\F}{{\mf F}}
\nc{\CA}{{\mathcal A}}
\nc{\CC}{{\mathcal C}}
\nc{\CE}{{\mathcal E}}
\nc{\CP}{{\mathcal P}}
\nc{\CR}{{\mathcal R}}
\nc{\CO}{{\mathcal O}}
\nc{\CQ}{{\mathcal Q}}
\nc{\CK}{{\mathcal K}}
\nc{\CZ}{{\mathcal Z}}
\nc{\bk}{{\mathbf k}}
\nc{\fg}{{\mathfrak g}}
\nc{\Ann}{\text{Ann}}
\nc{\Rad}{\text{Rad}}
\nc{\Res}{\text{Res}}
\nc{\Ind}{\text{Ind}}
\nc{\Ker}{\text{Ker}}
\nc{\id}{\text{id}}

\nc{\be}{\begin{equation}}
\nc{\ee}{\end{equation}}
\nc{\bee}{\begin{equation*}}
\nc{\eee}{\end{equation*}}
\nc{\Hom}{\text{Hom}}
\nc{\rh}{\vec}
\newenvironment{conds}{
                       
                        \begin{enumerate} }
                     {\end{enumerate} }

\newenvironment{num}{
                      
                      \begin{enumerate} }
                    {\end{enumerate} }

\newcommand{\Int}{{\mathbb Z}}
\newcommand{\Nat}{{\mathbb N}}

\newcommand{\real}{{\mathbb R}}
\newcommand{\mc}{\mathcal}

\newcommand{\ck}[1]{{#1}^\vee}

\begin{document}

\title{Root subsystems of loop extensions}
\author{M.~J.  Dyer and G.~I. Lehrer}
\address{Department of Mathematics 
\\ 255 Hurley Building\\ University of Notre Dame \\
Notre Dame, Indiana 46556, U.S.A.}
\email{dyer.1@nd.edu}
\address{School of Mathematics and Statistics\\
University of Sydney. Sydney. 2006. Australia.}
\email{gustav.lehrer@sydney.edu.au}
\begin{abstract}
 We completely classify the real  root subsystems of root systems of loop algebras of Kac-Moody Lie algebras.
This classification involves new notions of ``admissible subgroups'' of the coweight lattice of 
a root system $\Psi$, and ``scaling functions'' on $\Psi$. Our results generalise and simplify earlier
work on subsystems of  real affine root systems.
 \end{abstract}
\maketitle

\section{Introduction}\label{intro}

\subsection{Purpose of this work} Starting with the set $\Phi$ of real roots of a Kac-Moody Lie
algebra $\fg$, we define its loop extension $\wh\Phi$, which is a root system in the sense of being closed
under all reflections corresponding to its elements. The set $\wh\Phi$ is realised as a set
of roots of an extended loop algebra $\wh\fg$ of $\fg$ (see below).
The purpose of this work is to determine all root subsystems
of $\wh\Phi$. These are subsets of $\wh \Phi$ with the above closure property.
In the case where $\Phi$ is a finite crystallographic
root system, this amounts to the determination of the subsystems of the corresponding affine root system.

The subsystems are characterised in terms of certain explicit maps $\CZ:\Phi\to\CP(\Z)$, where $\CP(\Z)$
is the power set of $\Z$, whose image is a ``compatible'' collection of cosets in $\Z$, and whose support is a subsystem
$\Psi$ of $\Phi$. 
Our main result implies that such a function $\CZ$ with support $\Psi$
corresponds to a single coset in the coweight lattice $\Omega(\Psi)$ (Definition \ref{def:cwtlat})
of an ``admissible subgroup'' (Definition \ref{def:adsub}). These admissible subgroups turn
out to correspond precisely to ``scaling functions'' $M:\Psi\to\N$; generically, these are functions such that 
$\Psi_M:=\{M(\alpha)\inv\alpha\mid\alpha\in\Psi\}$ and $\ck\Psi_M:=\{M(\alpha)\ck\alpha\mid\alpha\in\Psi\}$
is a dual pair of crystallographic root systems. Such functions are determined numerically in \S\ref{s4}.

Our treatment here is a generalisation of \cite[Theorem 4]{DyL} to arbitrary root systems $\Phi$, and 
provides a more natural and more general
definition of ``admissible coweight lattices'' (cf. \S\ref{ss:last}) than that in \cite{DyL}.

\subsection{Background}
Let $\mathfrak  g$ be a Kac-Moody Lie algebra over $\C$ (see \cite{Kac} or \cite{Mo}).
We denote by $\wh {\mathfrak  g}$ the extension of the loop algebra
$\mathfrak {g}\otimes_{\C}\C[t,t^{-1}]$ by the derivation $d=\text{\rm Id}_{\mathfrak{g}}\otimes t\frac{d}{dt}$.
As $\C$-vector space,
we have 
\bee \wh{\mathfrak {g}}=\bigl(\mathfrak {g}\otimes _{\C}\C[t,t^{-1}] \bigr)\oplus \C d.  \eee
The Lie bracket on $\wh{\mathfrak g}$  is given by 
\bee [t^{m}x+\lambda d,t^{n}y+\mu d]=t^{m+n}[x,y]_{0}+n \lambda t^{n}y-m\mu t^{m}x \eee
for $x,y\in \mathfrak {g}$, $m,n\in \Int$ and $\lambda,\mu\in \C$, where $[-,-]_{0}$ is the bracket 
on $\mathfrak {g}$ and we write $t^{m}x=x\otimes t^{m}$, etc.

Let $\mathfrak{h}$ denote the standard Cartan subalgebra of $\mathfrak g$
and $\mathfrak{g}=\mathfrak{h}\oplus\bigoplus_{\alpha\in \Delta}\mathfrak{g}_{\alpha}$ be the root space 
decomposition  of $\mathfrak{g}$, where $\Delta\subseteq \mathfrak{h}^{*}$ is the root system and 
\bee \mathfrak{g}_{\alpha}=\mset{x\in \mathfrak{g}\mid [h,x]_{0}=\alpha(h)x\text{ \rm for all $h\in \mathfrak{h}$}}\eee 
for any $\alpha$ in the dual space $\mathfrak{h}^{*}$.
We define the abelian Lie subalgebra  $\wh{\mathfrak{h}}=\mathfrak{h}\oplus \C d$ of $\wh{\mathfrak{g}}$ and identify 
the dual space of $\wh{\mathfrak{h} }$ as $\wh{\mathfrak{h}}^{*}={\mathfrak{h}}^{*}\oplus \C \delta$ where 
$\delta(d)=1$, $\delta({\mathfrak{h}})=0$, $\mathfrak{h}^{*}(d)=0$.
Let $\wh \Delta:=\bigl((\Delta\cup\set{0}) +\Int\delta\bigr)\setminus \set{0}\subseteq \wh{\mathfrak{h}}^{*}$.
Then there is again a root space decomposition
 $ \wh {\mathfrak{g}}=\wh{\mathfrak{h}}\oplus\bigoplus_{\alpha\in \wh \Delta}\wh{\mathfrak{g}}_{\alpha}$ 
 where for any  $\alpha\in \Delta\cup\set{0}$ and $m\in \Int$ with $\alpha+m\delta\neq 0$, we have  
 $\wh{\mathfrak{g}}_{\alpha+m\delta}=t^{m}{\mathfrak{g}}_{\alpha}$.
 We call $\wh \Delta$ the loop root system of $\Delta$.
 
 Let $\Phi:=\Delta^{\text{\rm re}}\subseteq \Delta$ be the subset of real roots of $\Delta$
 and $\wh \Phi :=\Phi+\Int\delta\subseteq \wh \Delta$. We shall refer to
 the elements of $\wh \Phi$ as the real roots of $\wh{\mathfrak{g}}$.
 However, note that since $\wh\fg$ is not generally a Kac-Moody Lie algebra, this
 terminology is not standard.
 
 Define reflections 
 $s_{\beta}\colon \wh{\mathfrak{h}}^{*}\rightarrow \wh{\mathfrak{h}}^{*}$ for $\beta\in \wh \Phi$ as follows.
 Write $\beta=\alpha+m\delta$ where 
  $\alpha\in \Phi$ and $m\in \Int$, and let  
  $s_{\beta}=s_{\alpha+m\delta}\colon \phi \mapsto \phi -\phi(\ck \alpha)(\alpha+m\delta) $,
  where $\ck\alpha\in\mathfrak{h}$ is the coroot corresponding to $\alpha$.
  For $\beta\in \Phi$, $s_{\beta}$ restricts to the usual reflection in $\beta$ on $\mathfrak{h}^{*}$.
  We let $W:=\mpair{s_{\alpha}\mid \alpha\in \Phi}$ denote the Weyl group of $\mathfrak{g}$
  (a Coxeter group) and $\wh W:=\mpair{s_{\alpha}\mid \alpha\in \wh\Phi}$.
  This is not known to be  a Coxeter group in general. 
  We have $W\Phi=\Phi$, $W\Delta=\Delta$, $\wh W\wh \Phi=\wh \Phi$,
and   $\wh W\wh \Delta=\wh \Delta$.
  
  A subset $\Psi$ of $\Phi$ (resp., $\wh \Phi$) is called a root subsystem of $\Phi$
  (resp., $\wh \Phi$) if $s_{\alpha}(\beta)\in \Psi$ for all $\alpha,\beta\in \Psi$.
    It is easy to see that   the root subsystems of  $\wh \Phi$ are 
  of the form $\wh \Psi_{\mc{Z}}:=\mset{\alpha+m\delta\mid \alpha\in \Psi,m\in Z_{\alpha}}$ 
  where $\Psi$ is a root subsystem of $\Phi$, $\set{Z_{\alpha}}_{\alpha\in \Psi} $ is a 
  family of subsets (in fact, cosets of subgroups) of $\Int$, and $\CZ:\Psi\to\CP(\Z)$
  is the function which takes $\alpha\in\Psi$ to $Z_\alpha\subseteq\Z$.
  There are easily described compatibility conditions among the subsets $Z_{\alpha}$  
  of $\Int$ which ensure that $\Psi_{\mc{Z}}$ is a root subsystem of $\wh \Phi$.
  In this work, we show how these conditions imply that the subsets are defined in terms
  of combinatorial data involving the coweight lattice of $\Psi$, and its admissible subgroups
  (cf. \S\ref{s2}). We show in turn that the latter are described in terms of scaling functions
  on $\Psi$.
  
  In this way, we obtain 
  as our main result a 
  bijective  parameterisation   (Theorem \ref{thm:bij}) of root subsystems of $\wh \Phi$ by  the 
  root subsystems of $\Phi$ together with additional simple  combinatorial data.
  
  If $\mathfrak{g}$ is a semisimple complex Lie algebra, then $\Phi$ is a finite crystallographic
  (reduced) root system and $W$ is a finite  Weyl group.
  In this case, $\wh{\mathfrak g}$ is the quotient  (by the centre) of the untwisted affine Kac-Moody Lie algebra
  corresponding to $\mathfrak g$, 
  $\wh \Phi$ is the corresponding  affine root system, and $\wh W$ is the affine  Weyl group  of $\Phi$.
  We recover in this case a particularly simple and illuminating  proof of  the 
  parameterisation \cite[Theorem 4]{DyL} of the reflection subgroups  of the affine Weyl group $\wh W$
  (or equivalently, of root subsystems of  the real affine root system $\wh \Phi$).
  
 The main result is formulated  purely in terms of root systems, without any mention of Lie algebras, 
 and its proof requires only general facts about crystallographic  Coxeter groups and their root systems  
 which can be found in \cite{Bour}, \cite{Mo}, \cite{DyRef} and \cite{Semi}. 
 In  \S \ref{s1}, starting with the notion
 of a ``crystallographic based root datum'', we develop 
 the theory of loop extensions of root systems and define the
 notion of lifting. This is used in \S\ref{s2} to define the ``coweight lattice'' $\Omega(\Psi)$
 of a root system $\Psi$. In that section, we also define admissible subgroups  of $\Omega(\Psi)$
 and prove our main result, Theorem \ref{thm:bij}. 
 In \S\ref{s4} we define scaling functions of $\Psi$  and  explicitly determine them   when $\Psi$ has a simple
 system whose components are finite, but whose Weyl group may be infinite.
 Finally, \S\ref{s3} elucidates the application to root subsystems of affine Weyl groups,
 and explains how the present results relate to those of \cite{DyL}.

 \section{Loop root systems} \label{s1}
  \subsection{}  
Recall that a subset $\Gamma$ of a $\real$-vector space $L$ is said to be  positively independent
   if $\sum_{\alpha\in \Gamma}c_{\alpha}\alpha=0$ with $c_{\alpha}\in \real_{\geq 0}$ and all 
 \label{ss1} 
   but finitely 
  many $c_{\alpha}=0$ implies that all $c_{\alpha}=0$. We shall also require the notion of a generalised Cartan
  matrix, which is a possibly infinite matrix $C:=(c_{\alpha,\beta})_{\alpha,\beta\in \Pi}$ satisfying
  $c_{\alpha,\alpha}=2$,  and for $\alpha\neq \beta$, $c_{\alpha,\beta}\in \Int_{\leq 0}$
  with $c_{\alpha,\beta}=0$ iff $c_{\beta,\alpha}=0$.
 \begin{definition}\label{def:cbrd}
  A {\em crystallographic  based root datum} is a tuple 
  \bee B=(V,{\ck V}, \Pi,\ck \Pi,\mpair{-,-},\iota)\eee  satisfying the following conditions: 
  \begin{conds}\item $V$ and ${\ck V}$ are $\real$-vector spaces and $\mpair{-,-}\colon V\times {\ck V}\rightarrow \real$ 
  is a $\real$-bilinear form.
  \item $\Pi\subseteq V$ and $\ck \Pi\subseteq {\ck V}$ are positively independent subsets.
  \item $\iota\colon \Pi\rightarrow \ck \Pi$ is a bijection (denoted $\alpha\mapsto \ck \alpha$) such that the matrix
  $C:=(c_{\alpha,\beta})_{\alpha,\beta\in \Pi}$ where $c_{\alpha,\beta}:=\mpair{\alpha,\ck \beta}$,  is
  a generalised Cartan matrix (see above).
  \end{conds}
  \end{definition}
 The definition above is modelled on that of a based root datum in \cite{DyRig};  a similar but more 
 technical  notion of crystallographic root datum can be found in \cite[Chapter 5]{Mo} for example. 
 We often omit the word crystallographic in this paper since we will not consider more general based root data. 
  
  For $\alpha\in \Pi$, define $s_{\alpha}\colon V\rightarrow V$ by $v\mapsto v-\mpair{v,\ck\alpha}\alpha$.
  Let $W=\mpair{s_{\alpha}\mid \alpha\in \Phi}$, $\Phi:=W\Pi$ and define the set of {\em positive roots}
  by $\Phi_{+}:=\Phi\cap \real_{\geq 0}\Pi$.
  Similarly, we define the group $\ck W$ as the group generated by reflections acting  on ${\ck V}$,  
  and the positive roots $\ck \Phi_{+}\subseteq \ck \Phi\subseteq {\ck V}$.
  It is easily seen that there is a canonical isomorphism $w\mapsto \ck w$ from $W$ to $\ck W$, and that the map
  $\iota$ extends uniquely to a bijection $\Phi\rightarrow \ck \Phi$, which transforms the action of
  $w\in W$ to that of $\ck w\in \ck W$. This extension is denoted $\iota\colon \alpha\mapsto \ck \alpha$.
  
  The  matrix $C$ is called the Cartan matrix of  $B$ and  $W$ is the Weyl group of $B$. 
  Note that $W$ is a Coxeter group, with simple reflections $S:=\mset{s_{\alpha}\mid \alpha\in \Pi}$;
  the  Coxeter matrix of $(W,S)$ is determined by the Cartan matrix $C$ in a well known way 
  (see \cite[Prop. 3.13]{Kac} or \cite[5.1, Proposition 11]{Mo}).  The elements of $\Phi$, $\Phi_{+}$, $\Pi$ 
  respectively are called  roots, positive roots and simple roots of $B$.
  
 The map $\alpha\mapsto s_{\alpha}$ induces a bijection $\Phi_{+}\rightarrow T$ 
 where $T=\mset{wsw^{-1}\mid w\in W,s\in S}$ is the set of reflections of $(W,S)$. 
 Subgroups of $W$ generated by subsets of $T$ are called reflection subgroups of $(W,S)$.

  \subsection{} Let  $\wh V=V\oplus \real \delta$, $\ck {\wh V}={\ck V}\oplus \real \delta'$ and extend 
  $\mpair{-,-}$ to a bilinear form $\wh V\times\ck {\wh V}\rightarrow \real$ with
  $\mpair{\wh V,\delta'}=\mpair{\delta,\ck{\wh V}}=\set{0}$ and $\mpair{\delta,\delta'}=1$.
  
  \begin{definition}\label{def:phihat}
  The loop extension $\wh\Phi$ of $\Phi$ is defined as 
  $\wh \Phi=\Phi+\Int \delta$.
  \end{definition}
  For  any $\alpha+m\delta\in \wh \Phi$, define $\wh s_{\alpha}\colon \wh V \rightarrow \wh V$ by 
  $v\mapsto v-\mpair{v,\ck \alpha}(\alpha+m\delta)$.
  For any $\gamma\in {\ck V}$, define $t_{\gamma}\colon \wh V\rightarrow \wh V$ by
  $t_{\gamma}(v)=v+\mpair{v,\gamma}\delta$. It is easily checked that
  $\wh s_{\alpha+m\delta}=\wh s_{\alpha}t_{-m\ck \alpha}$.

For any subgroup $Q$ of ${\ck V}$, let $t_{Q}:=\mset{t_{v}\mid v\in Q}$; this is a group
of  linear transformations of $\wh V$, isomorphic to a quotient of $Q$.
The group generated by $\{\wh s_{\alpha}\mid\alpha\in \wh \Phi\}$
will be denoted $\wh W$. The group $W$ is identified as a subgroup of $\wh W$ by identifying $s_{\alpha}\in W$ with 
$\wh s_{\alpha}\in \wh W$, and it is straightforward that
$\wh W$ is the semidirect product $\wh W=W\ltimes t_{Q}$ 
where $Q=\Int \ck \Phi$. If $W$ is  a finite Weyl group, then $\wh W$ is an affine Weyl group; in general,
$\wh  W$ is not known to have the structure of a Coxeter group.
 In view of this identification, we henceforth write $s_{\alpha}$ instead of $\wh s_{\alpha}$ 
for all $\alpha\in \Psi$.

\subsection{}\label{ss:subsyst} For any subset $\Gamma$ of $\Phi$, 
define the reflection subgroup $W_{\Gamma}:=\mpair{s_{\alpha}\mid \alpha\in \Gamma}$ of $W$.

\begin{definition}\label{def:subsystem} A root subsystem $\Psi$  of $\Phi$ (resp., $\wh \Phi$) is a subset 
of $\Phi$ (resp., $\wh \Phi$) such that $\alpha,\beta\in \Psi$  implies $s_{\alpha}(\beta) \in \Psi$.
\end{definition}
There is a bijection $\Psi\mapsto W_{\Psi}:=\mpair{s_{\alpha}\mid \alpha\in \Psi}$ 
from the set of root subsystems of $\Phi$ to the set of reflection subgroups of $W$. 
It is not known if an analogous statement holds for $\wh \Phi$ and $\wh W$.

\subsection{}\label{ss:refsubgp} In the statement below we record certain general 
properties  of  root subsystems and their corresponding reflection groups. These will be crucial 
below.

\begin{proposition}\label{prop:subprops}
Let $B$ be a based root datum as in \text{\rm Definition \ref{def:cbrd}}, with corresponding root system $\Phi$ 
and positive system $\Phi^+$.
\begin{num}

\item Every root subsystem  $\Psi$ of $\Phi$ has a 
canonical simple system $\Gamma$, defined 
by $\Gamma\subseteq \Phi_{+}$,  $W_{\Gamma}=W_{\Psi}$ and $\mpair{\alpha,\ck \beta}\leq 0$ 
for all distinct $\alpha,\beta\in \Gamma$.

\item The pair $(W_{\Psi},\mset{s_{\alpha}\mid \alpha\in \Gamma})$ is a Coxeter system. 

\item There is a  based root datum  $B_{\Psi}:=(V,\ck V, \Gamma, \ck \Gamma, \mpair{-,-},\iota')$
with Weyl group $W_{\Psi}$ where $\iota'$ is the map $\alpha\mapsto \ck \alpha\colon \Gamma \rightarrow \ck \Gamma$. 
\end{num}
\end{proposition}

A proof of the analogous facts for a slightly different class of root systems may be found in \cite{Mo}.
Similar arguments may be applied here using the ideas in  Lemma \ref{ss:basis} below.
The results are also the  specialisation to crystallographic root systems  of facts in \cite{DyRig}
(see also  \cite{DyRef}).

\subsection{} The following simple  fact will be useful.
\begin{lemma} \label{lem:triv} Let  $\Psi$ be a root subsystem of $\Phi$ 
with canonical simple system $\Gamma$ (see \text{\rm Prop. \ref{prop:subprops}}).
Then  $\Psi$ is characterised as the unique  subset $\Psi'$ of $\Phi$ with 
$\Gamma\subseteq \Psi'\subseteq \Psi$ such that  $\alpha,\beta\in \Psi'$ implies 
$s_{\alpha}(\beta)\in \Psi'$.\end{lemma} 

\begin{proof} We have  $\Psi=W_{\Gamma}\Gamma=W_{\Psi}\Psi$.
Hence if $\Psi'$ is as stated, then $\Psi'=W_{\Psi'}\Psi'$ and $\Gamma\subseteq \Psi'\subseteq \Psi$.
Hence $W_\Gamma\Gamma(=\Psi)\subseteq W_{\Psi'}\Psi'(=\Psi')\subseteq W_\Psi\Psi=\Psi$.
 \end{proof}

\subsection{Lifting} \label{ss:basis} In general the canonical simple system $\Gamma$ of a based root datum may be 
linearly dependent. Such situations arise for subsystems, even when
the simple roots of the ambient system are linearly independent.
In the construction below, we introduce a lifting of a root subsystem to 
another one with the same generalised Cartan matrix, but where the simple roots are linearly independent.
This is essential for our classification of the subsystems.

 Fix a based root datum $B$ as in Definition \ref{ss1} with root system $\Phi$,
 and let $\Psi$ be any root subsystem of $\Phi$. 
  We associate to $\Psi$ a new based root datum $B^0_{\Psi}$ as follows.
  Let $\Gamma$ be the canonical simple system of $\Psi$ (Proposition \ref{prop:subprops}).
  
  Let $\Gamma_{0}$ and $\ck\Gamma_{0}$ be sets which are in canonical bijection with $\Gamma$
  and $\ck \Gamma$ respectively. Denote the bijections by $\alpha\in\Gamma\mapsto\alpha_0\in\Gamma_0$
  and $\ck\alpha\in\ck\Gamma\mapsto\ck\alpha_0\in\ck\Gamma_0$.
  Let $V_{0}$ and ${\ck V_0}$ be $\real$-vector spaces with bases $\Gamma_0$ and $\ck\Gamma_0$ respectively.
  Let $\iota_{0}\colon \Gamma_{0}\rightarrow \ck \Gamma_{0}$ be the bijection induced by $\iota$. 
  There is a unique bilinear form $\mpair{-,-}_{0}\colon V_{0}\times \ck V_0\rightarrow \real$
  such that $\mpair{\alpha_{0},\ck\beta_{0}}_{0}=\mpair{\alpha,\ck\beta}$ for $\alpha,\beta\in \Gamma$.
  Now define
  \[ B^0_{\Psi}:=(V_{0},{\ck V_0},\Gamma_{0},\ck \Gamma_{0}, \mpair{-,-}_{0},\iota_{0}).\]
  This is a based root datum, as is verified by checking the definition. We denote the root system
  of $B^0_{\Psi}$ by $\Psi_{0}$, its dual root system by $\ck \Psi_{0}$  and its Weyl group by $W^0$.
  Denote also by $\iota_{0}$ the unique $W^0$-equivariant extension of $\iota_0:\Gamma_{0}\rightarrow \ck \Gamma_{0}$
  to $\Psi_{0}$. This is a bijection $\Psi_{0}\rightarrow \ck \Psi_{0}$.
  There are $\real$-linear maps $L\colon V_{0}\rightarrow V$ and 
  $\ck L\colon {\ck V_0}\rightarrow {\ck V}$  such that $L(\alpha_{0})=\alpha$ and $\ck L(\ck \alpha_{0})=\ck \alpha$
  for all $\alpha\in \Gamma$. For $\alpha_0\in \Psi_0$, denote the reflection in $\alpha_0$ by $s^0_{\alpha_0}$.
  \begin{lemma}\label{lem:basis}
  \begin{num}
  \item  The map $s^0_{\alpha_{0}}\mapsto s_{\alpha}$ for $\alpha\in \Gamma$ extends to a 
  group isomorphism $W^0\rightarrow W_{\Psi}$. We therefore identify these two groups.
  \item $\mpair{L(v'),\ck L(u')}=\mpair{v',u'}_{0}$ for $v'\in V_{0}$ and $u'\in {\ck V_0}$.
  \item $wL(v')=L(wv')$ and $w\ck L(u')=\ck L(wu')$ for $w\in W_{\Psi}$, $v'\in V_{0}$ and $u'\in {\ck V_0}$
  \item $L$ and $\ck L$ restrict to bijections $\Psi_{0}\rightarrow \Psi$ and $\ck \Psi_{0}\rightarrow \ck{\Psi}$ 
  satisfying $\ck L(\iota_{0}\alpha)=\iota{L(\alpha)}$. Equivalently, for $\alpha\in\Psi$,
  $\ck L(\ck\alpha_0)=\ck\alpha$.
    \end{num}\end{lemma}
    \begin{remark} In view of the properties, we   denote the inverse bijections in 
    (c)  as $\alpha\mapsto \alpha_{0}\colon  \Psi\mapsto \Psi_{0}$ and 
    $\ck \alpha \mapsto \ck\alpha_{0}\colon \ck \Psi\mapsto \ck \Psi_{0}$ respectively, 
    $s^0_{\alpha}$ just as $s_{\alpha}$,  and the bijection $\Psi_{0}\rightarrow \ck\Psi_{0}$ as 
    $\alpha\mapsto \ck \alpha$, with no risk of confusion.\end{remark}
  \begin{proof} We provide only a sketch of the proof, as the argument is well known;  see e.g. 
  \cite[5.1, Proposition 1 and Lemma 2]{Mo} for details of a similar argument. 
  Part (a) is proved using the well known formula for the Coxeter matrix in terms of the Cartan 
  matrix, and the fact that the Coxeter matrix determines $W_\Psi$.
  Part (b) is trivial, and part (c) is proved by induction on the length of $w$. Part (d)  
  follows using (a)--(c) and the last paragraph of \ref{ss1}.
   \end{proof}

\section{Subsystems of loop root systems}\label{s2}

   \subsection{}  For subsets $A$, $B$ of any abelian group, define 
   $A+ B:=\mset{a+ b\mid a\in A,b\in B}$,  $A-B:=\mset{a-b\mid a\in A,b\in B}$ and $nA=\mset{na\mid a\in A}$ for $n\in \Int$.
   Note that $2A\neq A+A$ in general, but $A-nB=A+(-n)B$, ($n\in \Z$).  
   
      Any subset $\Sigma$  of $\wh\Phi$ (see Definition \ref{def:phihat}) can be written in the form 
      $\Sigma:=\mset{\alpha+n\delta\mid {\alpha\in \Phi}, n\in Z_{\alpha}}$ for unique subsets 
      $Z_{\alpha}$ of $\Int$, for $\alpha\in \Phi$.
   The next Lemma gives equations amongst the $Z_{\alpha}$ which are necessary and 
   sufficient for the corresponding set $\Sigma$ to be a root subsystem of $\wh \Phi$. 
 
 \begin{lemma} \label{lem:Z}(cf. \cite[Lemma 9]{DyL})   
 Let $\mc{Z}:\Phi\to\CP(\Z),\;\alpha\mapsto Z_{\alpha}$ be a function. 
 Then $\Sigma:=\mset{\alpha+n\delta\mid {\alpha\in \Phi}, n\in Z_{\alpha}}$ is a root subsystem of $\wh\Phi$
 if and only if for all $\alpha,\beta\in \Phi$, we have   \begin{equation}\label{Z} 
Z_{\beta}-\mpair{\beta,\ck\alpha}Z_{\alpha}\subseteq 
 Z_{s_{\alpha}(\beta)}.\end{equation}
\end{lemma}
 \begin{proof} 
 We have $s_{\alpha+m\delta}(\beta+n\delta)=s_{\alpha}(\beta)+(n-m\mpair{\beta,\ck\alpha})\delta$.
 The result  follows directly from  this and the definitions.
   \end{proof}
   
 \begin{definition}\label{def:supp}
 \begin{enumerate}
 \item A function $\mc{Z}:\Phi\to\CP(\Z),\;\alpha\mapsto Z_{\alpha}$ which satisfies \eqref{Z}
 will be called a {\em \rfp}
\item The {\em support} of a \rf  $\CZ$ is $\{\alpha\in\Phi\mid\CZ_\alpha\neq\emptyset\}$.
 \end{enumerate}
 \end{definition}
 
 \begin{lemma}\label{lem:sup}
 The support of any \rf  $\CZ:\Phi\to\CP(\Z)$ is a root subsystem of $\Phi$.
 \end{lemma}
 The proof is easy.
 
 \subsection{}   The root subsystems $\Sigma$ of $\wh\Phi$ may be  described in terms of solutions 
 of \eqref{Z} as follows. 
 
 \begin{corollary} \label{cor:Zsyst}   The root subsystems $\Sigma$
 of  $\wh \Phi$ are in bijective correspondence with {\rfs}\  $\mc{Z}:\Phi\to\CP(\Z)$.
    The correspondence attaches to a \rf  $\CZ$ with support $\Psi$ 
   the root subsystem  $\wh \Psi_{\mc{Z}}:=\mset{\alpha+n\delta\mid \alpha\in \Psi,
   n\in Z_{\alpha}}$ of  $\wh \Phi$. \end{corollary}
 
 \begin{proof} This is clear from Lemmas \ref{lem:Z} and \ref{lem:sup}.
 \end{proof}

 \subsection{}
 Our goal is to give explicit descriptions of all \rfs\ 
 with  support equal to a fixed subsystem $\Psi$ of $\Phi$.
 We therefore fix such a subsystem $\Psi$ until further notice.

  Let  $\Gamma$ be  the canonical simple system of $\Psi$.
  Consider the associated based root datum $B^0_{\Psi}$ as in Lemma \ref{lem:basis}, 
   with root system $\Psi_{0}$ which is naturally in bijection with $\Psi$ via the map 
   $\alpha\mapsto \alpha_{0}\colon \Psi\mapsto \Psi_{0}\subseteq V_{0}$ and with 
   $\Gamma_{0}=\mset{\alpha_{0}\mid \alpha\in \Gamma}$ as $\real$-basis. 

\begin{definition}\label{def:cwtlat} Let $\Int\Psi_{0}$ be the additive subgroup of $V_{0}$ generated by 
$\Psi_{0}$. This is a free abelian group with basis $\Gamma_{0}$. 
\begin{num}
\item The group $\Int\Psi_{0}$ is the {\em root lattice} of $\Psi_{0}$.
\item The {\em coweight lattice} of $\Psi$ is the abelian group
$\Omega(\Psi):=\Hom_{\Int}(\Int \Psi_{0},\Int)$. It is free abelian if $\Gamma$ is finite.
\end{num}
\end{definition}

   Evidently $\Omega({\Psi})$ is an abelian group naturally isomorphic to the group 
   $\Int^{\Gamma}$ of all functions $\Gamma\rightarrow \Int$ with pointwise operations, 
   under the correspondence $f\mapsto\bigl( \gamma\mapsto f(\gamma_{0})\colon \Gamma\rightarrow \Int\bigr)$ 
   for $f\in \Omega(\Psi)$. A similar definition attaches to $f\in \Omega(\Psi)$ a function
   $\Psi\to \Z$; note that this function need not respect the linear relations among
   the elements of $\Psi$. This construction provides a source of \rfs as follows.
   
   For any subset $X$ of $\Omega(\Psi)$ and any $\beta\in \Int\Psi_{0}$, write
   $X(\beta):=\mset{f(\beta)\mid f\in X}\subseteq\Z$. We shall study conditions which ensure that the map
   $\beta\mapsto X(\beta_{0})$ ($\beta\in \Psi$) defines a \rf  with support $\Psi$.
   
   There is a natural $W_{\Psi}$-action on $\Omega(\Psi)$ 
   defined by $(wf)(\alpha)=f(w^{-1}\alpha)$ for all $f\in \Omega(\Psi)$, $w\in W_{\Psi}$ and $\alpha\in \Int \Psi_{0}$.
   
   Our first Lemma shows that the set of \rfs $\mc{Z}:\Phi\to\CP(\Z)$, 
   with support $\Psi$ in which 
   each $Z_{\alpha}=\CZ(\alpha)$ for $\alpha\in \Psi$ is a singleton subset $\set{p_{\alpha}}$ of  $\Int$ 
   is in natural bijection with $\Omega(\Psi)$.

\begin{lemma} \label{lem:p} A  function  $p:\Psi\to\Z$ 
given by $\alpha\mapsto p_{\alpha}$ satisfies 
\be\label{eq:relsing}
p_{\beta}-\mpair{\beta,\ck\alpha}p_{\alpha}=p_{s_{\alpha}(\beta)} 
\ee
for  all $\alpha,\beta\in \Psi$
if and only if there is some $f\in \Omega(\Psi)$ such that $p_{\alpha}=f(\alpha_{0})$ for all $\alpha\in \Psi$.
In that case, $f$ and $p$ are uniquely determined by the restriction $p_{\vert \Gamma }\colon  \Gamma\to \Z$,
which may be arbitrary. The function $p$ may be thought of as a \rf with support $\Psi$, whose value at
$\alpha\in\Psi$ is $\{p_\alpha\}$.
\end{lemma}
\begin{proof} Fix a function $p\mapsto p_\alpha:\Psi\to \Z$.
There is a unique  $f\in \Omega(\Psi)$ such that $p_{\alpha}=f(\alpha_{0})$ for all 
$\alpha\in \Gamma$.  Let $p':\Psi\to\Z$ be the function corresponding to $f$ as above,
i.e., $p'(\alpha)=p'_{\alpha}=f(\alpha_{0})$ for all $\alpha\in \Psi$.
Since $(s_{\alpha}(\beta))_{0}=\beta_{0}-\mpair{\beta,\ck \alpha}\alpha_{0}$
and $f$ is $\Z$-linear on $\Psi_0$, we have 
$p'_{\beta}-\mpair{\beta,\ck\alpha}p'_{\alpha}=p'_{s_{\alpha}(\beta)} $ for  all $\alpha,\beta\in \Psi$. 

We show that if $p$ satisfies the relation \eqref{eq:relsing}, then $p=p'$.
Suppose $p_{\beta}-\mpair{\beta,\ck\alpha}p_{\alpha}=p_{s_{\alpha}(\beta)} $ for  all $\alpha,\beta\in \Psi$.
Let  $p''_{\beta}:=p_{\beta}-p'_{\beta}$ for $\beta\in \Psi$.  Then
$p''_{\beta}-\mpair{\beta,\ck\alpha}p''_{\alpha}=p''_{s_{\alpha}(\beta)} $ for  all $\alpha,\beta\in \Psi$ 
and $p''_{\beta}=0$ for all $\beta\in \Gamma$. 
Note that if $\alpha,\beta\in \Psi$ with $p''_{\alpha}=p''_{\beta}=0$, then $p''_{s_{\alpha}(\beta)}=0$. 
It follows by Lemma \ref{lem:triv}  that $p''_{\beta}=0$ for all $\beta\in \Psi$, as required.

The above argument also shows that if $p$ arises as above from an arbitrary element $f\in\Omega(\Psi)$,
then $p$ satisfies  \eqref{eq:relsing}. Uniqueness follows also from the argument above.
\end{proof}

\subsection{} Next we study  the \rfs \;\;$\mc{Z}:\alpha\mapsto Z_\alpha$ with support $\Psi$ in which
each $Z_{\alpha}$ is a subgroup, denoted $n_{\alpha}\Int$, of $\Int$, where $n_{\alpha}\in \Nat$.

\begin{lemma} \label{lem:n} Let $N:\Psi\to\N$ be a function, and write $N(\alpha)=n_\alpha$.
Then the following conditions on the integers $n_\alpha$ are equivalent:
\begin{conds} \item For all $\alpha,\beta\in \Psi$,
  $n_{\beta}\Int-\mpair{\beta,\ck\alpha}n_{\alpha}\Int\subseteq n_{s_{\alpha}(\beta)}\Int $
  \item    $n_{\beta}\vert \mpair{\beta,\ck\alpha}n_{\alpha}$ for all $\alpha,\beta\in \Psi$,
and $n_{\alpha}=n_{\beta}$ for all $\alpha,\beta\in \Psi$ with $\beta\in W_{\Psi}\alpha$.
 \item  For  all $\alpha,\beta\in \Psi$,
$n_{\beta}\Int-\mpair{\beta,\ck\alpha}n_{\alpha}\Int= n_{s_{\alpha}(\beta)}\Int $. 
\end{conds}
\end{lemma}
\begin{proof} Suppose that (i) holds. The condition is equivalent to the conditions:
for all $\alpha,\beta\in \Psi$, we have $n_{s_{\alpha}(\beta)}\vert \mpair{\beta,\ck \alpha}\alpha$
and $ n_{s_{\alpha}(\beta)}\vert n_{\beta}$. This implies by symmetry that 
$n_{\beta}= n_{s_{\alpha}(\beta)}$ for $\alpha,\beta\in \Psi$. It follows that
$n_{\alpha}=n_{\beta}$ whenever $\alpha,\beta\in \Psi$ are in the same $W_{\Psi}$-orbit.
Consequently (ii) holds. Clearly, (ii) implies (iii) and (iii) implies (i).
\end{proof}
 
 \begin{lemma}\label{lem:admsubgp}  Let $M:\Gamma\to\N$ be a function, and write $M(\alpha)=m_\alpha$.
 Suppose that for all $\alpha,\beta\in \Gamma$, $m_{\beta}\vert \mpair{\beta,\ck\alpha}m_{\alpha}$.
 Define
 \[X_{M}:=\mset{f\in  \Omega(\Psi)\mid f(\alpha_{0})\in m_{\alpha}\Int
 \text{ \rm  for all $\alpha\in \Gamma$}}.
 \] 
 This is the largest 
 subgroup $X$ of $\Omega(\Psi)$ such that
 $X(\alpha_{0})\subseteq m_{\alpha}\Int $ for all $\alpha\in \Gamma$. 
 Define the function $N:\Psi\to\N$ ($\alpha\mapsto n_\alpha$) by
 $n_{\alpha}\Int=X_{M}(\alpha_{0})$. Then
 \begin{num}\item   The subgroup $X_{M}$ of $\Omega(\Psi)$ is $W_{\Psi}$-stable.  
  \item  The function $N$ is the unique extension of $M$ to a 
  function $\Psi\to\N$
  satisfying  the equivalent conditions of $\text{\rm Lemma \ref{lem:n}}$.
 \end{num}
 \end{lemma}
 \begin{proof} To prove (a), it will suffice to show that if $f\in X_{M}$ and $\alpha\in \Gamma$ 
 then $s_{\alpha}f\in X_{M}$ i.e. that $(s_{\alpha}f)(\beta_{0})\in \Int m_{\beta}$ for 
 all $\beta\in \Gamma$. But 
 \bee (s_{\alpha}f)(\beta_{0})=f((s_{\alpha}(\beta))_{0})
 =f(\beta_{0})-\mpair{\beta,\ck \alpha}f(\alpha_{0})\in 
 m_{\beta}\Z+\mpair{\beta,\ck \alpha}m_{\alpha}\Z=m_{\beta}\Z,
 \eee 
 as required. This proves (a). 
 
 For any $\alpha\in \Gamma$,  $X_{M}$ contains  $m_{\alpha}\omega_{\alpha}$ where 
 the ``fundamental coweight'' $\omega_{\alpha}\colon 
 \Int \Psi_{0}\rightarrow \Int$ is the $\Int$-linear map determined
 by $\omega_{\alpha}(\beta_{0})=m_{\alpha}\delta_{\alpha,\beta}$ for
 all $\beta\in \Gamma$. Hence $m_{\alpha}\in X_{M}(\alpha_{0})\subseteq m_\alpha\Int$ 
 and  therefore $X_{M}(\alpha_{0})=m_{\alpha}\Int$. It follows that
 $n_{\alpha}=m_{\alpha}$ for $\alpha\in \Gamma$ 
 i.e. $N$ is an extension of $M$. That $N$ satisfies 
 the conditions of Lemma \ref{lem:n} follows from the first part of
 Corollary \ref{cor:star} (see Remark \ref{rem:lemma7}). 
 Finally, the uniqueness of the extension follows from the second  condition in 
 Lemma \ref{lem:n}(ii).
  \end{proof}
  \begin{remark} The above proof shows that  a function $M\colon \alpha\mapsto m_{\alpha}$
  satisfying the conditions in Lemma \ref{lem:admsubgp}  has the additional property that 
  $m_{\alpha}=m_{\beta}$ whenever $\alpha,\beta\in \Gamma$ are in the same $W_{\Psi}$-orbit. 
  This may  be seen more directly using a well-known criterion \cite[Ch IV, \S 1, Proposition 3]{Bour} for simple roots 
  to be in the same $W_{\Psi}$-orbit and the formula for the Coxeter matrix of $W_{\Psi}$
  in terms of the Cartan matrix for $B^{0}_{\Psi}$.
  \end{remark}
 
  \begin{corollary}\label{cor:star}
  Let $N:\Psi\to\N$ be a function, and write $N(\alpha)=n_\alpha$.
   The conditions \text{\rm (i)--(iii)} of \text{\rm Lemma \ref{lem:n}} on the integers $n_\alpha$
   are equivalent to the following statement:
   $(*)$ There is a $W_{\Psi}$-stable subgroup  $X$ of $\Omega(\Psi)$ such that 
$n_{\alpha}\Int=X(\alpha_{0})$ for all $\alpha\in \Psi$.
  \end{corollary}
  \begin{proof}
  Suppose that $X$ is any $W$-stable subgroup of $\Omega(\Psi)$; define integers
  $n_{\alpha}\in \Nat$  for $\alpha\in \Psi$, by $X(\alpha_{0})=n_{\alpha }\Int $.
  For any $w\in W_{\Psi}$ and $\alpha\in \Psi$,  $X((w\alpha)_{0})=(w^{-1}X)(\alpha_{0})=X(\alpha_{0})$,
  which shows that $n_{\alpha}=n_{w\alpha}$.
  Moreover for $\alpha,\beta\in \Psi$ and $f\in X$, we have
  $f((s_{\alpha}(\beta))_{0})=f(\beta_{0})-\mpair{\beta,\ck \alpha}f(\alpha_{0})$.
  Fix $\alpha\in\Psi$ and choose $f$ so $f(\alpha_{0})=n_{\alpha}$. 
  Then 
  \[\mpair{\beta,\ck \alpha}n_{\alpha}=
  f(\beta_{0})-f((s_{\alpha}(\beta))_{0})\in n_{\beta}\Int +n_{s_{\alpha}(\beta)}\Int=
  n_{\beta}\Int\] 
  and so $n_{\beta}\vert \mpair{\beta,\ck\alpha}n_{\alpha}$.
  This shows that the integers $n_\alpha$ satisfy (ii) of Lemma \ref{lem:n}. 
  \begin{remark}\label{rem:lemma7}
  The above implication is used in the proof of Lemma \ref{lem:admsubgp}. Note that
  its proof does not involve Lemma \ref{lem:admsubgp}.
  \end{remark}

 Conversely, suppose (ii) of Lemma \ref{lem:n} holds for the integers $n_\alpha$.
 Applying Lemma \ref{lem:admsubgp}  to the function $M$
   with $m_{\alpha}=n_{\alpha}$ for all $\alpha\in \Gamma$, we see that
   $n_{\alpha}=X_{M}(\alpha_{0})$ for all $\alpha\in \Psi$ where 
   $X_{M}$ is a $W_{\Psi}$-stable subgroup of $\Omega(\Psi)$.
  \end{proof}
  
   \begin{definition}\label{def:adsub} We call a subgroup 
    $X$ of the coweight lattice $\Omega(\Psi)$ an {\em admissible 
    subgroup} if there is a function $N:\Psi\to\N$ satisfying the conditions
    of Lemma \ref{lem:n} such that 
    \[X=\mset{f\in  \Omega(\Psi)\mid f(\alpha_{0})\in n_{\alpha}\Int
 \text{ \rm  for all $\alpha\in \Psi$}}.
 \] 
 \end{definition}
 
 \begin{remark}
     Note that for any admissible subgroup $X$, if $M:\Gamma\to \N$ is the restriction of 
     the function $N$ then $X$ is determined by $M$, since $X=X_M$
    as in Lemma \ref{lem:admsubgp}.
 \end{remark}

  The following statement is a summary of the results of this subsection.
  
  \begin{proposition}\label{prop:sum}
  Let $\CZ:\Phi\to\CP(\Z)$ be a \rf with support $\Psi$ such that for each $\alpha\in\Psi$,
  $\CZ(\alpha)$ is a subgroup of $\Z$. Then there is a unique admissible subgroup $X\subseteq\Omega(\Psi)$
  such that for all $\alpha\in\Psi$, $\CZ(\alpha)=X(\alpha_0):=\{x(\alpha_0)\mid x\in X\}$.   
  \end{proposition}

  \subsection{} \label{ss:Zsubgp}  In order to utilise Lemma \ref{lem:p} and
  Proposition \ref{prop:sum}, we prove the following.   
   
   \begin{lemma} \label{lem:Zsubgp} Let  $\mc{Z}\colon \Phi\to \CP(\Z)$
   be a \rf with support $\Psi$ and let $Z_{\alpha}:=\CZ(\alpha)$.
      Suppose that $0\in {Z}_{\alpha}$ for all $\alpha\in \Gamma$.
   Then $Z_{\alpha}$ is a subgroup of $\Int$ for all $\alpha\in \Psi$.
   \end{lemma}
   \begin{proof} Note that if $\alpha,\beta\in \Psi$ and $0\in Z_{\alpha}\cap Z_{\beta}$, 
   then by \eqref{Z}, $0\in Z_{s_{\alpha}(\beta)}$. It follows from Lemma \ref{lem:triv} 
   that if $0\in {Z}_{\alpha}$ for all $\alpha\in \Gamma$, then
   $0\in Z_{\alpha}$ for all $\alpha\in \Psi$. Also $Z_{\alpha}-2Z_{\alpha}\subseteq Z_{-\alpha}$.
The left hand side contains $-Z_{\alpha}$, so we get that $-Z_{\alpha}\subseteq Z_{-\alpha}$.
Replacing $\alpha$ by $-\alpha$, we see that $Z_{-\alpha}=-Z_{\alpha}$
and substituting this into $Z_{\alpha}-2Z_{\alpha}\subseteq Z_{-\alpha}$, we see that 
$2Z_{\alpha}-Z_{\alpha}= Z_{\alpha}$.
From this and $0\in Z_{\alpha}$, it follows 
that  $Z_{\alpha}=-Z_{\alpha}$, and 
\be\label{eq:Z2}
2Z_{\alpha}+Z_{\alpha}= Z_{\alpha}.
\ee
Then by induction on $n$ we see that 
$nZ_{\alpha}\subseteq Z_{\alpha}$ for all $n\in \Int$.  

If $Z_{\alpha}=\set{0}$, 
then we are done. Otherwise, let
$n_{\alpha}$ be the least positive element of $Z_{\alpha}$.  
We have $n_{\alpha}\Int \subseteq Z_{\alpha}$, and we claim equality holds. 
Take any element $n\in Z_{\alpha}$.
Then $n=q2n_\alpha+r$ for unique integers $q,r$ with $-n_\alpha< r\leq n_\alpha$. 
From \eqref{eq:Z2}, it follows that $r\in Z_{\alpha}$, whence $-r\in Z_{\alpha}$.
By the minimal nature of $n_{\alpha}$, we conclude that $r\in \set{0,n_\alpha}$,
and hence that $n\in n_\alpha\Z$.

Hence  $Z_{\alpha}=n_{\alpha}\Int$ which is a subgroup of $\Int$ as claimed.   
   \end{proof}

\subsection{The main theorem} We start with a key definition.
\begin{definition}\label{def:R}
Let $\mc{R}$ denote the set of all 
pairs $(\Psi,Y)$ where $\Psi$ is a root subsystem 
of $\Phi$ and $Y\in \Omega(\Psi)/X$ is a coset of some admissible subgroup
$X$ of  $\Omega(\Psi)$.
\end{definition}
\begin{theorem}   \label{thm:bij} There is a bijection between $\mc{R}$ and the 
set of all \rfs which to any element $(\Psi,Y)\in \mc{R}$, attaches  the unique \rf
$\CZ:\Phi\to\CP(\Z)$ with support $\Psi$ such that for $\alpha\in\Psi$,
$\CZ(\alpha)=Z_\alpha=Y(\alpha_0)$
 \end{theorem}
 \begin{remark}
 In view of Corollary \ref{cor:Zsyst} this shows that $\CR$ is in natural bijection with 
 the root subsystems of $\wh \Phi$.
 \end{remark}
\begin{proof} 

We first show that the function $\CZ$ is a \rf.
Let $X$ be an admissible subgroup of $\Omega(\Psi)$, and $x\in \Omega(\Psi)$ such that $Y=x+X$.
Let  $Z'_{\alpha}:=\set{x(\alpha_{0})}$ and $Z''(\alpha):=X(\alpha_{0})$ for $\alpha\in \Psi$.
By Lemma \ref{lem:p} and Corollary \ref{cor:star}, $Z'_{\alpha}$ and $Z''_{\alpha}$ each satisfy \eqref{Z} for all  
$\alpha,\beta\in \Psi$. Hence  $Z_{\alpha}:=Y(\alpha_0)=Z'_{\alpha}+Z''_{\alpha}$ satisfies \eqref{Z} for all
$\alpha,\beta\in \Psi$. 

Next, we show that any \rf $\CZ$ is realised in this way. Let $\Psi$ be the support of $\CZ$.
Let $\Gamma$ be the canonical simple system of $\Psi$.  Choose $p_{\alpha}\in Z_{\alpha}$ for all 
$\alpha\in \Gamma$, and extend to a \rf $p$ with support $\Psi$ as in
Lemma \ref{lem:p}; as in {\it loc.~cit.}, we have $p_{\alpha}=x(\alpha_{0})$ (for all $\alpha\in\Psi$)
for some $x\in \Omega(\Psi)$. Set 
$Z'_{\alpha}:=Z_{\alpha}-p_{\alpha}$ for all $\alpha\in \Psi$.
Then $\mc{Z}':\alpha\mapsto\mset{Z'_{\alpha}}\;\;({ \alpha\in \Psi})$ satisfies \eqref{Z} for all
$\alpha,\beta\in \Psi$, and so defines a \rf with support $\Psi$, such that
$0\in Z'_{\alpha}$ for all $\alpha\in \Gamma$. Hence $0\in Z'_{\alpha}$ for all $\alpha\in \Psi$ by Lemma \ref{lem:Zsubgp}. 
Again by  Lemma \ref{lem:Zsubgp} and Proposition \ref{prop:sum},
there is an admissible subgroup $X$ of $\Omega(\Psi)$ such that $Z'_{\alpha}=X(\alpha_{0})$ 
for all $\alpha\in \Psi$. Hence for all 
$\alpha\in \Psi$, $Z_{\alpha}=p_{\alpha}+Z'_{\alpha}=x(\alpha_{0})+X(\alpha_{0})=Y(\alpha_{0})$ 
where $Y:=x+X\in \Omega(\Psi)/X$. 

The theorem will now follow from the following assertion: if  $(\Psi_{i},Y_{i})\in \mc{R}$ ($i=1,2$) 
have equal root functions $\CZ$ as constructed above,
then $(\Psi_{1},Y_{1})=(\Psi_{2},Y_{2})$.
To prove this, note first that equal root functions have the same support,
so $\Psi:=\Psi_{1}=\Psi_{2}$. Now write $Y_{i}=x_{i}+X_{i}$, where for $i=1,2$,
$x_{i}\in \Omega(\Psi)$ and $X_{i}$ is an admissible subgroup of $\Omega(\Psi)$.
Let $\alpha\in \Psi$.
By assumption, we have $Y_{1}(\alpha_{0})=Y_{2}(\alpha_{0})$
i.e. $Z_{\alpha}:=x_{1}(\alpha_{0})+X_{1}(\alpha_{0})=x_{2}(\alpha_{0})+X_{2}(\alpha_{0})$.
Hence $Z_{\alpha}-Z_{\alpha}=X_{1}(\alpha_{0})=X_{2}(\alpha_{0})$.
By Definition \ref{def:adsub}, we have 
\bee X_{i}=
\mset{x\in \Omega(\Psi)\mid  x(\alpha_{0})\in X_{i}(\alpha_{0})\text{ \rm for all $\alpha\in \Psi$}}
\eee  
which is independent of $i$.
Hence $X:=X_{1} =X_{2}$. Now let $x:=x_{1}-x_{2}\in \Omega(\Psi)$. Since
$x(\alpha_{0})\in X(\alpha_{0})$ for all $\alpha\in \Psi$,  it follows, again by Definition \ref{def:adsub}, 
that $x\in X$ i.e. $Y_{1}=x_{1}+X_{1}=x_{2}+X_{2}=Y_2$ as required. This completes the proof.
\end{proof}
\begin{remark} It is easy to see from the proof of Theorem \ref{thm:bij}  that 
any \rf $\CZ\colon \alpha\mapsto Z_{\alpha}$   satisfies  the following stronger version of \eqref{Z}:
    \begin{equation}\label{Z=} 
Z_{\beta}-\mpair{\beta,\ck\alpha}Z_{\alpha}= 
 Z_{s_{\alpha}(\beta)}, \qquad \alpha, \beta\in \Phi.\end{equation}
Note that \eqref{Z=} also  holds as an equation in $V$ if
 $Z_{\beta }$ is replaced by $\beta\in V$ for all $\beta\in \Phi$.
 The equation \eqref{Z=} can be interpreted  with the $Z_{\alpha}$ elements of an additive  monoid equipped 
 with a suitable  operation by $\Int$; the elements $Z_{\beta}$ of the monoid may be thought of 
 as  ``abstract roots''.    The  proof  in this section of the Theorem implicitly describes the 
 solutions of \eqref{Z=}
 in various such  monoids: the monoids  of singleton subsets of   $\Int$, subgroups of $\Int$ and 
 subsets of $\Int$. Similarly, it can be solved with the $Z_{\alpha}$ discrete subsets of $\real$. 
 Similar equations exist for non-crystallographic root systems.\end{remark}
 
 \section{Explicit description of subsystems; scaling functions.}\label{s4}
 
 We have seen (Theorem \ref{thm:bij}) that the root subsystems of $\wh\Phi$ are in bijection with
 the set $\CR$ of pairs $(\Psi,Y)$, where $\Psi$ is a subsystem of $\Phi$ and $Y$ is a coset
 $x+X$ of some admissible subgroup $X$ of $\Omega(\Psi)$. Taking the classification of the subsystems
 $\Psi$ of $\Phi$ as known, the crucial element of the classification is therefore that of the
 admissible subgroups of $\Omega(\Psi)$. In this section we provide such a classification
 when each component $\Gamma_i$ of $\Gamma$ is finite. This means simply that 
 $|\Gamma_i|<\infty$; it may still happen that 
 $\Gamma_i$ is not of finite type, i.e. that the root system $W_{\Gamma_i}\Gamma_i$ is infinite. 
 
 \subsection{Scaling functions and scaled root systems}  \label{ss:Mfamily}
  Let  $\Psi$ be a root subsystem of $\Phi$ with canonical simple system  
  $\Gamma$ and associated root datum 
  $$
  B_{\Psi}:=(V,\ck V, \Gamma, \ck \Gamma, \mpair{-,-},\iota')
  $$ 
  as in \ref{ss:refsubgp}.
     By Lemma \ref{lem:admsubgp}, the admissible subgroups of $\Omega(\Psi)$ are in 
     natural bijective correspondence with certain functions  $M:\Gamma\to \N$.
 
 \begin{definition}\label{def:scalf} With the above notation, a {\em scaling function} 
 is a function $M:\Gamma\to\N$ ($\alpha\mapsto m_\alpha$)
 which satisfies the requirement that for all $\alpha,\beta\in \Gamma$, 
 $m_{\beta}\vert\mpair{\beta,\ck \alpha}m_{\alpha}$. 
 \end{definition}
 
 Each scaling function on $\Gamma$ extends 
 uniquely to a function $N:\Psi\to\N$ $(\alpha\mapsto n_\alpha)$, which satisfies the conditions 
 in Lemma \ref{lem:n}(ii). This extension will be referred to as a scaling function on $\Psi$.
 
  We begin with some remarks concerning scaling functions. First observe that $M$ is a scaling
  function for $\Gamma$ if and only if its restriction to each of its components is a scaling
  function for that component. We therefore henceforth confine our attention to the case
  where $\Gamma$ is irreducible. In that case, if  $m_{\alpha}=0$ for some $\alpha\in \Gamma$, 
 then $m_{\alpha}=0$ for all $\alpha\in \Gamma$. Hence we also assume that no value of $M$ is $0$.
 
If $M:\Gamma\to\Nat,\;\;\alpha\mapsto m_\alpha$ is any function with 
$m_{\alpha}\neq 0$ for all $\alpha\in\Gamma$, then $M$ is a scaling function
if and only if $\mpair{m_{\alpha}^{-1}\alpha, m_{\beta}\ck\beta}\in \Z$
for all $\alpha,\beta\in \Gamma$. In that case, there 
 is an associated  (crystallographic)  based root datum
  $B_{M}:=(V,\ck V, \Delta, \ck \Delta, \mpair{-,-},\iota'')$ defined as follows: $\Delta:=\set{m_{\alpha}^{-1}\alpha
\mid \alpha\in \Gamma}$ and $\ck \Delta:=\set{m_{\alpha}\ck\alpha
\mid \alpha\in \Gamma}$ and $\iota''(m_{\alpha}^{-1}\alpha)=m_{\alpha}  \ck \alpha$ for all
$\alpha\in \Gamma$. The Weyl group of $B_{M}$ is equal to $W_{\Psi}$. The root system of $B_{M}$ is
$\Psi_{M}:=\mset{n_{\alpha}^{-1}\alpha\mid \alpha\in \Psi}$, where $\alpha\mapsto n_\alpha$ is the
unique extension of $M$ to $\Psi$ discussed above. The terminology ``scaling function''
is suggested by the fact that this root system
is  obtained by rescaling the roots in  $\Psi$. The admissible subgroup 
$X_{M}$ of $\Omega(\Psi)$ corresponding to $M$ is then naturally isomorphic to 
the coweight lattice $\Omega(\Psi_{M})$.

   \subsection{Determination of scaling functions} We next determine the non-zero scaling functions 
   $M\colon \alpha\mapsto m_{\alpha}$ in the case that $\Gamma$ is finite and irreducible.
  Note that one always has the ``trivial'' constant scaling functions, and that
  if $\alpha\mapsto m_\alpha$ is a scaling function and $k\in\N$ then $\alpha\mapsto km_\alpha$
  is again a scaling function.
  
  For any positive prime $p\in \Int$, let $\nu_{p}\colon \mathbb{Q}\rightarrow \mathbb{Z}\cup \set{\infty}$
  be the $p$-adic valuation; then $\nu_{p}(0)=\infty$ and $\nu_{p}(x)=n\in \Int$
  if $x=p^{n}\frac{a}{b}$ where $a,b\in \Int$, $n\in \Int$ and $p\nmid a$, $p\nmid b$.
     The condition that $M$ be a scaling function is equivalent to   
     $\nu_{p}(m_\alpha)\in \Nat$ for all  $\alpha\in \Gamma$ and 
  \be \label{ineq3} \nu_{p}(m_{\beta})\leq \nu_{p}(\mpair{\beta,\ck \alpha})+\nu_{p}(m_{\alpha})\ee
  for all $\alpha,\beta\in \Gamma$, and all primes $p$.  
  It is convenient to first consider solutions of the inequalities \eqref{ineq3}
  with $\nu_{p}(m_\alpha)\in \Int$ (corresponding to $m_{\alpha}\in \mathbb{Q}^{*}=\mathbb{Q}\setminus\set{0}$).
  Given any $\alpha,\beta\in \Gamma$, by the irreducibility of $\Gamma$ there is  a sequence
  $\alpha=\alpha_{0},\alpha_{1},\ldots, \alpha_{n}=\beta$ in $\Gamma$ such 
  that $n\in \Nat$  and  $\mpair{\alpha_{i+1},\ck\alpha_{i}}\neq 0$ for all $i=1,\ldots, n$.
  The inequalities \eqref{ineq3} imply that  
  \be\label{ineq2} \nu_{p}(m_{\beta})\leq\sum_{i=1}^{n} \nu_{p}(\mpair{\alpha_{i+1},\ck 
  \alpha_{i}})+\nu_{p}(m_{\alpha})\ee for  any such sequence.
  We see from this that there exist integers $c_{\alpha,\beta,p}\in \Nat$
  such that  the inequalities \eqref{ineq3} hold  if and only if
   \be \label{ineq}\nu_{p}(m_{\beta})\leq \nu_{p}(m_{\alpha})+c_{\alpha,\beta,p}\ee 
   for all
   $\alpha,\beta\in \Gamma$ and all primes $p$.
 The integer $c_{\alpha,\beta,p}$ is the minimum 
 value of $\sum_{i=1}^{n} \nu_{p}(\mpair{\alpha_{i+1},\ck 
     \alpha_{i}})$ over all  sequences $\alpha=\alpha_{0},\alpha_{1},\ldots \alpha_{n}=\beta$ as above.
 Note that if, for all $\alpha,\beta\in \Gamma$, $p$  does not divide $\mpair{\alpha,\ck \beta}\neq 0$, then one has 
 $c_{\alpha,\beta,p}=0$ for all $\alpha,\beta\in \Gamma$ and so $\nu_{p}(\alpha)$ is constant for $\alpha\in \Gamma$.
 Hence, to classify the non-constant scaling functions, it suffices to check \eqref{ineq} for a finite set
 $\CQ$ of primes, which is determined by $\Gamma$.

     For fixed $p$, regard \eqref{ineq}  as a set of inequalities for
   $\nu_{p}(m_{\alpha})\in \Int$ for all $\alpha\in \Gamma$. 
    Fix any element $\gamma\in \Gamma$.
    For any  solution  $\nu_{p}(m_{\alpha})=d_{\alpha,p}$ of \eqref{ineq}, 
    we have another solution $\nu_{p}(m_{\alpha}):=d_{\alpha,p}-d_{\gamma,p}$,
    and $\nu_{p}(m_{\gamma})=0$.
    But there are only finitely 
    many solutions $\nu_{p}(m_{\alpha})=d_{\alpha,p}$  of \eqref{ineq} with 
    $\nu_p(m_{\gamma})=0$,  and they may be effectively determined,  since the inequalities
   \eqref{ineq} imply that the set $\{\vert \nu_{p}(m_\alpha)- \nu_{p}(m_{\gamma})\vert 
   \mid \alpha\in\Gamma\}$ of absolute values is  bounded. 
     Write these finitely many solutions as 
     $\nu_{p}(m_{\alpha})=d_{\alpha,p}^{(k)}(\in\Z)$ for $\alpha\in \Gamma$,
      for $k=1,\ldots, N_{p}$. From each such solution, one obtains a solution of 
      \eqref{ineq} with all $\nu_{p}(m_{\alpha})\in \Nat$, and at least one
        $\nu_{p}(m_{\alpha})=0$, given by
          $\nu_{p}(m_{\alpha})=e_{\alpha,p}^{(k)}$ where
           $e_{\alpha,p}^{(k)}:=d_{\alpha,p}^{(k)}-d^{(k)}_{p}$ and 
            $d_{p}^{(k)}:=\min\set{d_{\alpha,p}^{(k)}\mid \alpha\in \Gamma}$.
            It is clear from the above that the solutions of 
          \eqref{ineq} with all $\nu_{p}(m_{\alpha})\in \Nat$ are precisely the 
          solutions of the form  $\nu_{p}(m_{\alpha})=e_{\alpha,p}^{(k)}+K_p$
           for all $\alpha\in \Gamma$, where $K_p\in \Nat$ is any
           non-negative integer, independent of $\alpha$, and  
           $ k=1,\ldots, N_{p}$.

   Allowing $p$ to vary over the finite set $\CQ$ it is now clear that we
   have proved the following statement.
   
   \begin{proposition}\label{prop:scalf}
   Let $\Gamma\subseteq\Phi$ be a finite, irreducible simple system, and retain the above notation.
   For any sequence $\bk=(k_p)_{p\in\CQ}$ with $1\leq k_p\leq N_p$, define 
   \[
   m_\alpha^{(\bk)}=\prod_{p\in\CQ}p^{e_{\alpha,p}^{(k_{p})}}.
   \]
   
   The scaling functions $M:\Gamma\to\N$ are precisely the functions of the form
   $M(\alpha)=nm_\alpha^{(\bk)}$, where $\bk$ is a sequence as above, and $n\in\N$.
 \end{proposition}
   
   This follows from the above discussion, from which it is clear that
   the non-zero scaling functions $\alpha\mapsto m_{\alpha}$   are given by
   \be 
   m_{\alpha}=\prod_{p\in\CQ} p^{e_{\alpha,p}^{(k_{p})}+K_p}\prod_{p\not\in\CQ}p^{K_p},
   \ee  
   for $\alpha\in \Gamma$, with $k_{p}=1,\ldots, N_{p}$, $K_p\in \Nat$ and all but finitely many 
    $K_p=0$.

   In this way, we see that when $\Gamma$ is finite, there are finitely many ``basic'' non-zero scaling
     functions (in fact, $\prod_{p}N_{p}$ of them), and any non-zero  scaling  
    function is  given by multiplying one of the basic ones by the positive integer
 $\prod_{p}p^{K_p}$.

 \begin{remark}  
 (1) The determination of all scaling functions in Proposition \ref{prop:scalf} applies whenever
all components $\Gamma_i$ of $\Gamma$ are finite. Note that this does not imply that 
$\Gamma_i$ is of finite type, i.e. $W_{\Gamma_i}$ may be infinite.
 
 (2)  If $\Phi$  has only finite and affine type components, then 
 the above techniques in conjunction with  results in \cite{DyL} (giving 
 descriptions of the canonical simple systems of  root subsystems of $\Phi$) 
 are enough to effectively determine the root subsystems of $\wh \Phi$.
 We describe the results concretely in the case in which $\Phi$ is of finite type in  \S5.    
 
 (3) In general, a parametrisation of simple systems  of the root subsystems of $\Phi$ 
 is not explicitly known. Further, such simple systems $\Gamma$ may have some infinite  rank 
 components  $\Gamma_{i}$  even if $\Phi$ is of finite rank. This complicates the explicit determination 
 of scaling functions in the general case.  
 \end{remark}

 \section{Root subsystems of affine root systems}
 \label{s3}  In this final section, we indicate  how the preceding results, may be applied
 in the case of root systems of finite Weyl groups, to give 
 an alternative proof of a parameterisation of root subsystems of real affine root systems
 described in \cite[\S 5]{DyL}.

\subsection{}   Henceforth, we assume that $\Phi$ is finite, and, for simplicity, that $\Phi$ spans 
$V$ as $\real$ vector space. Then $\Phi$ is a finite reduced  (crystallographic) root system with 
Weyl group $W$ in the sense of \cite{Bour}. Moreover, $\wh W$ is naturally isomorphic to the affine 
Weyl group $W^{a}$ of $W$, which may be realised as 
the  group of affine transformations of $V$ generated by $W$ and the group of translations by elements 
of the coroot lattice of $\Phi$ (see \cite{Bour}). The root system $\wh \Phi$ is 
the system of real roots of an untwisted affine  Kac-Moody Lie algebra
corresponding to $\Phi$ (\cite{Kac}) and $\wh W$ identifies with the (affine) Weyl group of this (affine) root system. 
In this case, root subsystems of $\wh \Phi$ correspond bijectively to  
reflection subgroups of $W^{a}$.
 
To make the description of root subsystems of $\wh \Phi$ explicit, we must  
describe explicitly all  the  scaling functions $M\colon \alpha\mapsto m_{\alpha}$ for $\Gamma$.
\subsection{}  
  Fix the subsystem $\Psi$ of $\Phi$ and let $\Gamma$ be its canonical simple system.   
  Let $\Gamma_{i}$  ($i=1,\ldots, n$) be the components of $\Gamma$. We let $k_{i}$ denote 
  the maximum ratio of squares of  root lengths of two roots of 
  $\Gamma_{i}$ (with  respect to any $W$-invariant  inner product on $V$).
  It is well known (see \cite[Ch VI, \S 1, no. 3]{Bour}) that $k_{i}\in \set{1,2,3}$;
  in fact, $k_{i}=3$ if $\Gamma_{i}$ is of type $G_{2}$,
  $k_{i}=2$ if $\Gamma_{i}$ is of type $F_{4}$, $B_{n}$ or $C_{n}$ with $n\geq 2$, 
  and otherwise $k_{i}=1$.  If there is only one root length ($k_i=1$), we say
  that all roots are both short and long.
          
        \begin{proposition} \label{Zcosetsize1}  A function
        $M\colon \Gamma\to \N, \;\;(\alpha\mapsto m_{\alpha})$  is a scaling function     
        if and only if  for each $i=1,\ldots, n$, at least one of the following three properties holds:
        \begin{conds}\item $m_{\alpha}=0$ for all $\alpha\in \Gamma_{i}$.
        \item there is some $q_{i}\in \Nat_{> 0}$ such that $m_{\alpha}=q_{i}$ for all $\alpha\in \Gamma_{i}$.
        \item there is some $q_{i}\in \Nat_{> 0}$ such that $m_{\alpha}=q_{i}$ for all 
        short $\alpha\in \Gamma_{i}$   and $m_{\alpha}=q_{i}k_{i}$ for all long $\alpha\in \Gamma_{i}$.     
        \end{conds}
               \end{proposition}
        \begin{remark} If $k_{i}=1$, (ii) and (iii) are equivalent.
        \end{remark}
        \begin{proof}
         It is known from \cite{Bour} that there are at most two root lengths in $\Gamma_{i}$, and 
         that the roots of each length form a single $W_{\Gamma_{i}}$-orbit.
        For each $i$, there are therefore  natural numbers  $q_{i}$  and $q_{i}'$
        (which are equal if $k_i=1$) such that $m_{\alpha}=q_{i}$ for all short $\alpha\in \Gamma_{i}$ and 
        $m_{\alpha}=q_{i}'$ for all long $\alpha\in \Gamma_{i}$.
         Using the well  known description \cite{Bour} of rank two root systems of $\Psi$, the condition 
        $m_{\alpha}\vert \mpair{\alpha,\ck\beta} m_{\beta}$ holds  for $\alpha,\beta\in \Psi$ if and
        only if  for each $i$,  $q_{i}\vert q_{i}'$ and $q_{i}'\vert k_{i }q_{i}$. Since $k_{i}\in \set{1,2,3}$,
        the result follows. \end{proof}
        \subsection{}\label{ss:last} In \cite[\S 5]{DyL}, there is attached to each root subsystem $\Psi$ of $\Phi$, 
        with canonical simple system $\Gamma$, a family of so-called ``admissible coweight lattices for 
        $\Psi$''. Using the  preceding Proposition, they are
         seen to be precisely the coweight lattices of the root systems $\Psi_{M}$ attached to  scaling 
         functions $M$ for $\Gamma$ in \S\ref{ss:Mfamily}, regarded as subsets of the coweight lattice of 
         $\Psi$ i.e. they are the admissible subgroups of $\Omega(\Psi)$ as defined in this  paper.
      Using this identification,   Theorem \ref{thm:bij}  specializes (in the case when $\Phi$ is finite) to the 
  parameterisation in \cite[Theorem 4]{DyL}  of root subsystems of affine root systems $\wh \Phi$.
  Several of the other results of \cite[\S 5]{DyL} for affine root systems extend to the loop root 
  systems $\wh \Phi$ in general.     

      \bibliography{loopsub}
\bibliographystyle{plain}

\nobreak

\end{document}